\documentstyle[12pt,amscd]{amsart}

\def\romm{\vrule height1em width0cm depth0.4em}
\def\rommm{\vrule height1.2em width0cm depth0.6em}

\let\:=\colon
\def\i{^{-1}}

\def\CC{{\bold C}}
\def\ZZ{{\bold Z}}
\def\QQ{{\bold Q}}
\def\PP{{\bold P}}

\def\I{{\cal I}}

\def\E{{\cal E}}

\def\C{{\cal C}}
\def\Z{{\cal Z}}

\def\OO{{\cal O}}

\def\L{{\cal L}}
\def\x{\times}
\def\*{\otimes}
\def\v{^{\vee}}
\def\iso{\simeq}

\def\sub{\subseteq}

\def\Hom{\operatorname{Hom}}
\def\Ext{\operatorname{Ext}}
\def\End{\operatorname{End}}
\def\Hilb{\operatorname{Hilb}}

\def\GL{\operatorname{GL}}

\def\rank{\operatorname{rank}}

\def\dsum{\operatornamewithlimits{\oplus}}
\def\ytre{\displaystyle{\operatornamewithlimits{\wedge}}}
\def\Ker{\operatorname{Ker}}
\def\Coker{\operatorname{Coker}}

\def\Supp{\operatorname{Supp}}
\def\codim{\operatorname{codim}}

\def\OP#1{{\OO_{\PP^#1}}}

\def\l{{\lambda}}

\newtheorem{thm}{Theorem}[section]

\newtheorem{lem}[thm]{Lemma}
\newtheorem{prop}[thm]{Proposition}

\theoremstyle{definition}
\newtheorem{defn}[thm]{Definition}

\theoremstyle{remark}
\newtheorem{rem}[thm]{Remark}
\newtheorem{notation}{Notation}
\newtheorem{ack}{Acknowledgements}

\numberwithin{equation}{section}

\newcommand{\thmref}[1]{theorem~\ref{#1}}
\newcommand{\propref}[1]{proposition~\ref{#1}}
\newcommand{\lemref}[1]{lemma~\ref{#1}}
\newcommand{\secref}[1]{\S\ref{#1}}

\begin{document}
\bibliographystyle{plain}

\title
{ Bott's formula and enumerative geometry}

\author{Geir Ellingsrud}
\address{Mathematical Institute\\University of Oslo\\P.~O.~Box~1053\\
         N--0316 Oslo, Norway}
\email{ellingsr@@math.uio.no}
\author{Stein Arild Str{\o}mme}
\address{Mathematical Institute\\University of Bergen\\All\'eg 55\\
         N--5007 Bergen, Norway}
\email{stromme@@mi.uib.no}
\thanks{Copyright \copyright 1994 G.~Ellingsrud and S. A. Str{\o}mme.
All rights reserved.}

\date{November 1994}
\dedicatory{Dedicated to the memory of Alf Bj{\o}rn Aure, 1955--1994}


\maketitle

\begin{abstract}
We outline a strategy for computing intersection numbers on smooth
varieties with torus actions using a residue formula of Bott.  As an
example, Gromov-Witten numbers of twisted cubic and elliptic quartic
curves on some general complete intersection in projective space are
computed.  The results are consistent with predictions made from mirror
symmetry computations.  We also compute degrees of some loci in the
linear system of plane curves of degrees less than 10, like those
corresponding to sums of powers of linear forms, and curves carrying
inscribed polygons.
\end{abstract}

\section{Introduction}
One way to approach enumerative problems is to find a suitable
complete parameter space for the objects that one wants to count, and
express the locus of objects satisfying given conditions as a certain
zero-cycle on the parameter space.  For this method to yield an
explicit numerical answer, one needs in particular to be able to
evaluate the degree of a given zerodimensional cycle class.  This is
possible in principle whenever the numerical intersection ring (cycles
modulo numerical equivalence) of the parameter space is known, say in
terms of generators and relations.

Many parameter spaces carry natural actions of algebraic tori, in
particular those coming from projective enumerative problems.  In 1967,
Bott \cite{Bott-1} gave a residue formula that allows one to express
the degree of certain zero-cycles on a smooth complete variety with an
action of an algebraic torus in terms of local contributions supported
on the components of the fixpoint set.  These components tend to have
much simpler structure than the whole space; indeed, in many interesting
cases, including all the examples of the present paper, the fixpoints
are actually isolated.

We show in this note how Bott's formula can be effectively used to
attack some enumerative problems, even in cases where the rational
cohomology ring structure of the parameter space is not known.

Our first set of applications is the computation of the numbers of
twisted cubic curves (theorems \ref{main} and \ref{main2}) and
elliptic quartic curves (\thmref{main3}) contained in a general
complete intersection and satisfying suitable Schubert conditions.
The parameter spaces in question are suitable components of the
Hilbert scheme parameterizing these curves.  These components are
smooth, by the work of Piene and Schlessinger \cite{Pien-Schl} in the
case of cubics, and Avritzer and Vainsencher \cite{Avri-Vain} in the
case of elliptic quartics.

The second set of applications is based on the Hilbert scheme of
zero-dimensional subschemes
of $\PP^2$, which again is smooth by Fogarty's work \cite{Foga-1}. These
applications deal with the degree of the variety of sums of powers of
linear forms in three variables (\thmref{main4}) and Darboux curves
(\thmref{main5}).

\begin{ack}
Part of this work was done at the Max-Planck-Institut f\"ur Mathematik
in Bonn during the authors' stay there in the spring of 1993. We would
like to thank the MPI for this possibility. We would also like to
express our thanks to S. Katz, J. Le Potier, D. Morrison, and A.
Tyurin for raising some of the problems treated here and for many
stimulating discussions.
\end{ack}

\subsection{Main results}
The first theorem deals with the number of twisted cubics on a general
Calabi-Yau threefold which is a complete intersection in some projective
space.  There are exactly five types of such threefolds: the quintic in
$\PP^4$, the complete intersections $(3,3)$ and $(2,4)$ in $\PP^5$, the
complete intersection $(2,2,3)$ in $\PP^6$ and finally $(2,2,2,2)$ in
$\PP^7$.

\begin{thm}\label{main}
For the general complete intersection Calabi-Yau threefolds, the numbers
of twisted cubic curves they contain are given by the following table:

\smallskip
\begin{center}
\begin{tabular}
{|l|c|c|c|c|c|}\hline
Type of complete intersection\romm
& $5$    & $4,2$     & $3,3$    & $3,2,2$  &$2,2,2,2$\\ \hline
Number of twisted cubics \romm
&$317206375$ &$15655168$ &$6424326$ &$1611504$ &$416256$\\ \hline
\end{tabular}
\end{center}
\end{thm}

In the case of a general quintic in $\PP^4$, the number of rational
curves of any degree was predicted by Candelas et al.\ in
\cite{Cand-Gree-Ossa-Park}, and the
cubic case was verified by the authors in \cite{Elli-Stro-3}.  In
\cite{Libg-Teit}
Libgober and Teitelbaum predicted the corresponding numbers for the
other Calabi-Yau complete intersections.  Our results are all in
correspondence with their predictions.

Greene, Morrison, and Plesser \cite{Gree-Morr-Ples} have also predicted
certain numbers of rational curves on higher dimensional Calabi Yau
hypersurfaces.  Katz \cite{Katz-2} has verified these numbers for lines
and conics for hypersurfaces of dimension up to 10.  The methods of the
present paper have allowed us to verify the following numbers.  All but
the last one, $N_3^{1,1,1,1}(8)$, have been confirmed by D.~Morrison
(privat communication) to be consistent with \cite{Gree-Morr-Ples}.

\begin{thm}\label{main2}
For a general hypersurface $W$ of degree $n+1$ in $\PP^n$ $(n\le8)$ and
for a partition $\lambda=(\lambda_1\ge\dots\ge\lambda_m>0)$ of $n-4$,
the number $N_3^{\lambda}(n)$ of twisted cubics on $W$ which meet $m$
general linear subspaces of codimensions $\lambda_1+1,\dots,\lambda_m+1$
respectively is given as follows:
\medskip
\begin{center}
\begin{tabular}{|l|l|l||l|l|l|}\hline
\rommm$n$ & $\lambda$ & $N_3^{\lambda}(n)$&
\romm$n$ & $\lambda$ & $N_3^{\lambda}(n)$ \\ \hline
\romm$4$&{}  &$317206375$ &
\romm$7$&$1,1,1$  & $12197109744970010814464$ \\ \hline
\romm$5$& $1$ & $6255156277440$ & $8$&
\romm$4$  & $897560654227562339370036$ \\ \hline
\romm$6$& $2$       & $30528671745480104$ &$8$&
\romm$3,1$    & $17873898563070361396216980$ \\ \hline
\romm$6$& $1,1$    & $222548537108926490$ &
\romm$8$&$2,2$    & $33815935806268253433549768$ \\ \hline
\romm$7$&$3$      & $154090254047541417984$ &
\romm$8$&$2,1,1$  & $174633921378662035929052320$ \\ \hline
\romm$7$&$2,1$    & $2000750410187341381632$ &
\romm$8$&$1,1,1,1$& $957208127608222375829677128$ \\ \hline
\end{tabular}
\end{center}
\end{thm}

The number $N_3^{1,1,1,1}(8)$ is not related to mirror symmetry as far
as we know; Greene et.al.~ get numbers only for partitions with at most 3
parts.  Our methods also yield other numbers not predicted (so far!) by
physics methods: for example, there are 1345851984605831119032336
twisted cubics contained in a general nonic hypersurface in $\PP^7$ (not
a Calabi-Yau manifold).

A similar method can be used to compute the number of elliptic
quartic curves on general Calabi-Yau complete intersections. Here are
the results for some hypersurfaces, which we state without proof:

\begin{thm}\label{main3}
The number of quartic curves of arithmetic
genus 1 on a general hypersurface of degree $n+1$ in $\PP^n$ are for
$4\le n\le13$ given by the following table. These curves are all smooth.
\smallskip
\begin{center}
\begin{tabular}{|l|l|}\hline
$n$\rommm & Smooth elliptic quartics on a general hypersurface of
                              degree $n+1$ in $\PP^n$\\ \hline
$4$&\romm   $3718024750$\\ \hline
$5$&\romm   $387176346729900$\\ \hline
$6$&\romm   $81545482364153841075$\\ \hline
$7$&\romm   $26070644171652863075560960$\\ \hline
$8$&\romm   $12578051423036414381787519707655$\\ \hline
$9$&\romm   $8760858604226734657834823089352310000$\\ \hline
$10$&\romm  $8562156492484448592316222733927180351143552$\\ \hline
$11$&\romm  $11447911791501360069250820471811603020708611018752$\\ \hline
$12$&\romm
$20498612221082029813903827233942127541022477928303274152$\\ \hline
$13$&\romm
$48249485834889092561505032612701767175955799366431126942036480$\\ \hline
\end{tabular}
\end{center}
\end{thm}

This computation uses the description given in \cite{Avri-Vain} of the
irreducible component of the Hilbert scheme of $\PP^3$ parameterizing
smooth elliptic quartics.  This Hilbert scheme component can be
constructed from the Grassmannian of pencils of quadrics by two explicit
blowups with smooth centers, and one may identify the fixpoints for the
natural action of a torus in a manner analogous to what we carry out for
twisted cubics in this paper.  For another related construction, see
\cite{Meur-1}, which treats curves in a weighted projective space.

The number of elliptic quartics on the general quintic threefold was
predicted by Bershadsky et.al.\ \cite{Bers-Ceco-Oogu-Vafa}.  Their
number, 3721431625, includes singular quartics of geometric genus 1.
These are all plane binodal quartics, and their number is
$1185*2875=3406875$ by \cite{Vain-1}.  Thus the count of
\cite{Bers-Ceco-Oogu-Vafa} is compatible with the number above.

Recently, Kontsevich \cite{Kont-1} has developed a technique for
computing numbers of rational curves of {\em any\/} degree, using the
stack of stable maps rather than the Hilbert scheme as a parameter
space.  He also uses Bott's formula, but things get more complicated
than in the present paper because of the presence of non-isolated
fixpoints in the stack of stable maps.

The next theorem deals with plane curves of degree $n$ whose equation
can be expressed as a sum of $r$ powers of linear forms.  Let $PS(r,n)$
be the corresponding subvariety of $\PP^{n(n+3)/2}$.  Then $PS(r,n)$ is
the $r$-th secant variety of the $n$-th Veronese imbedding of $\PP^2$.
Let $p(r,n)$ be the number of ways a form corresponding to a general
element of $PS(r,n)$ can be written as a sum of $r$ $n$-th powers if
this number is finite, and 0 otherwise.  The last case occurs if and
only if $\dim(PS(r,n))$ is less than the expected $3r-1$.  We don't know
of an example where $p(r,n)>1$ if $PS(r,n)$ is a proper subvariety.  If
$p(r,n)=1$, then $p(r,n')=1$ for all $n'\ge n$. It is easy to see that
$p(2,n)=1$ for $n\ge3$.

\begin{thm} \label{main4} Assume that $n\ge r-1$ and $2\le r\le8$.
Then $p(r,n)$ times the degree of $PS(r,n)$ is $s_r(n)$, where
{\allowdisplaybreaks
 \begin{align}
   2\,s_2(n)= &\,n^4-10\,n^2+15\,n-6, \notag \\
   3!\,s_3(n)= &\,n^6-30\,n^4+45\,n^3+206\,n^2-576\,n+384,\notag \\
   4!\,s_4(n)= &\,n^8-60\,n^6+90\,n^5+1160\,n^4-3204\,n^3-5349\,n^2+26586
              \,n-23760,\notag\\
   5!\,s_5(n)= &\,n^{10}-100\,n^8+150\,n^7+3680\,n^6-10260\,n^5-
               52985\,n^4+\notag\\*
               &\,224130\,n^3+127344\,n^2-1500480\,n + 1664640,\notag\\
   6!\,s_6(n)=
     &\,n^{12}-150\,n^{10}+225\,n^9+8890\,n^8-25020\,n^7-244995\,n^6
     +1013490\,n^5+\notag
     \\*&\,2681974\,n^4-17302635\,n^3+1583400\,n^2+101094660\,n
     -134190000,\notag\\
   7!\,s_7(n)=&\,n^{14}-210\,n^{12}+315\,n^{11}+
          18214\,n^{10}-51660\,n^9-802935\,n^8+\notag\\*
          &\,3318210\,n^7+17619994\,n^6-102712365\,n^5
          -136396680\,n^4+\notag\\*
          &\,1498337820\,n^3-872582544\,n^2-7941265920\,n
          +12360418560,\notag\\
   8!\,s_8(n)=&\,n^{16}-280\,n^{14}+420\,n^{13}+
              33376\,n^{12}-95256\,n^{11}
                 -2134846\,n^{10}+\notag\\*
          &\,8858220\,n^9+75709144\,n^8-
          427552020\,n^7-1332406600\,n^6+\notag\\*
          &\,11132416680\,n^5+5108998089\,n^4
          -145109970684\,n^3+\notag\\*
          &\,144763373916\,n^2 +713178632880\,n-1286736675840.\notag
 \end{align}
 }
\end{thm}

For example, $s_5(4)=0$; this corresponds to the classical but
non-obvious fact that not all ternary quartics are sums of five fourth
powers.  (Those who are are called Clebsch quartics; they form a
hypersurface of degree 36).

Note in particular that $s_3(3)=4$.  It is classically known that
$PS(3,3)$ is indeed a hypersurface of degree 4, its equation is the
so-called $S$-invariant \cite{Salm}.  It follows that $p(3,3)=1$, and
hence that $p(3,n)=1$ for $n\ge3$.

Only the first few of these polynomials are reducible:
$s_r(r-1)=0$ for $r\le5$, but
the higher $s_r$ in the table are irreducible over $\QQ$.

Note that the formulas of the theorem are not valid unless $n\ge r-1$.
For example, a general quintic is uniquely expressable as a sum of
seven fifth powers (cfr.~the references in \cite{Muka-1}), while
$s_7(5)$ is negative.

The final application quite similar.  A {\em Darboux curve\/} is a
plane curve of degree $n$ circumscribing a complete $(n+1)$-gon (this
terminology extends the one used in \cite{Bart-1}).  This means that
there are distinct lines $L_0,\dots,L_n$ such that $C$ contains all
intersection points $L_i\cap L_j$ for $i<j$.  Equivalently, there are
linear forms $\ell_0,\dots,\ell_n$ such that the curve is the divisor
of zeroes of the rational section $\sum_{i=0}^n \ell_i\i$ of
$\OP2(-1)$.  Let $D(n)$ be the closure in $\PP^{n(n+3)/2}$ of the
locus of Darboux curves.  Let $p(n)$ be the number of inscribed
$(n+1)$-gons in a general Darboux curve, if finite, and 0 otherwise.

\begin{thm} \label{main5}
For $n=5,6,7,8,9$, the product of $p(n)$ and the degree of the Darboux locus
$D(n)$ is $2540, 583020, 99951390, 16059395240, 2598958192572$,
respectively.
\end{thm}

We have no guess as to what $p(n)$ is; it might well be 1 for
$n\ge5$.  It is always positive for $n\ge5$ by an argument of Barth's
\cite{Bart-1}.
For $n\le4$ it is 0.  For $n\le3$, all curves are Darboux.  For
$n=4$, Darboux curves are L\"uroth quartics, and form a degree 54
hypersurface \cite{Morl,LePo,Tyur-1}.


\section{Bott's formula}\label{C}
Let $X$ be a smooth complete variety of dimension $n$, and assume that
there is given an algebraic action of the multiplicative group $\CC^*$
on $X$ such that the fixpoint set $F$ is finite.  Let $\E$ be an
equivariant vector bundle of rank $r$ over $X$, and let
$p(c_1,\dots,c_r)$ be a weighted homogenous polynomial of degree $n$
with rational coefficients, where the variable $c_i$ has degree $i$.
Bott's original formula \cite{Bott-1} expressed the degree of the
zero-cycle $p(c_1(\E),\dots,c_r(\E))\in H^{2n}(X,\QQ)$ purely in terms
of data given by the representations induced by $\E$ and the tangent
bundle $T_X$ in the fixpoints of the action.

Later, Atiyah and Bott \cite{Atiy-Bott-1} gave a more general formula,
in the language of equivariant cohomology.  Its usefulness in our
context is mainly that it allows the input of Chern classes of several
equivariant bundles at once.  Without going into the theory of
equivariant cohomology, we will give here an interpretation of the
formula which is essentially contained in the work of Carrell and
Lieberman
\cite{Carr-Lieb-1,Carr-Lieb-2}.

To explain this, first note that the $\CC^*$ action on $X$ induces, by
differentiation, a global vector field $\xi\in H^0(X,T_X)$, and
furthermore, the fixpoint set $F$ is exactly the zero locus of $\xi$.
Hence the Koszul complex on the map $\xi\v\:\Omega_X \to \OO_X$ is a
locally free resolution of $\OO_F$.  For $i\ge0$, denote by $B_i$ the
cokernel of the Koszul map $\Omega_X^{i+1}\to \Omega_X^{i}$.
It is well known that $H^j(X,\Omega_X^i)$ vanishes for $i\ne j$, see
e.g. \cite{Carr-Lieb-1}.
Hence there are natural exact sequences for all $i$:
$$
0 \to H^i(X,\Omega_X^i) @>p_i>> H^i(X,B_i) @>r_i>> H^{i+1}(X,B_{i+1}) \to0.
$$
In particular, there are natural maps
$q_i=r_{i-1}\circ\dots\circ r_0\:H^0(F,\OO_F)\to H^i(X,B_i)$.

\begin{defn} Let $f\: F \to \CC$ be a function and $c\in
H^i(X,\Omega_X^i)$ a non-zero cohomology class. We say that $f$ {\em
represents\/} $c$ if $q_i(f)=p_i(c)$.
\end{defn}

For each $i\ge-1$, put $A_i=\ker q_{i+1}$. Then
$$
0=A_{-1}\sub\CC=A_0\sub A_1\sub A_2\sub \dots \sub A_n = H^0(F,\OO_F)
$$
is a filtration by sub-vector spaces
of the ring of complex-valued functions on $F$.
The filtration has the property that
$A_i A_j \sub A_{i+j}$, and the associated graded ring
$\dsum A_{i}/A_{i-1}$ is naturally isomorphic to the cohomology ring
$H^*(X,\CC)\iso \dsum H^i(X,\Omega_X^i)$.
(In \cite{Carr-Lieb-2}, the filtration is constructed as coming from one of
the spectral sequences associated to hypercohomology of the Koszul
complex above.)

An interesting aspect of this is that cohomology classes can be
represented as functions on the fixpoint set.  The representation is
unique up to addition of functions coming from cohomology classes of
lower degree (i.e., lower codimension). Since the algebra of functions
on a finite set is rather straightforward, this gives an efficient way
to evaluate zero-cycles, provided that 1) we know how to describe a
function representing a given class, and 2) we have an explicit formula
for the composite linear map
$$
\epsilon_X\: H^0(\OO_F) @>q_n>> H^{n}(X,\Omega_X^n)
@>\text{res}_X>\iso> \CC.
$$
These issues are addressed in the theorem below.

Let $\E$ be an equivariant vector bundle of rank $r$ on $X$.  In each
fixpoint $x\in F$ the fiber of $\E$ splits as a direct sum of
one-dimensional representations of $\CC^*$; let
$\tau_1(\E,x),\dots,\tau_r(\E,x)$ denote the corresponding
weights, and for all integers $k\ge0$, let $\sigma_k(\E,x)\in\ZZ$ be the
$k$-th elementary symmetric function in the $\tau_i(\E,x)$.

\begin{thm}\label{bott1}
Let the notation and terminology be as above. Then
\begin{enumerate}
\item
The $k$-th Chern class $c_k(\E)\in H^k(X,\Omega_X^k)$ of $\E$ can be
represented by the function $x\mapsto \sigma_k(\E,x)$.
\item
For a function $f\: F \to \CC$, we have
$\epsilon_X(f) =
  \displaystyle{\sum_{x\in F} \frac{f(x)}{\sigma_n(T_X,x)}}.$
\end{enumerate}
\end{thm}

\begin{pf} See \cite[equation~3.8]{Atiy-Bott-1}, and \cite{Carr-Lieb-2}.
\end{pf}

Note that the function $\sigma_k(\E,-)$ depends on the choice of a
$\CC^*$-linearisation of the bundle $\E$, whereas the Chern class
$c_k(\E)$ it represents does not.

\section{Twisted cubics} \label{A}
Let $\Hilb^{3t+1}_{\PP^n}$ be the Hilbert scheme parameterizing
subschemes of $\PP^n$ ($n\ge3$) with Hilbert polynomial $3t+1$, and let
$H_n$ denote the irreducible component of $\Hilb^{3t+1}_{\PP^n}$
containing the twisted cubic curves.  Recall from \cite{Pien-Schl} that
$H_3$ is smooth and projective of dimension $12$.  Any curve
corresponding to a point of $H_n$ spans a unique 3-space, hence $H_n$
admits a fibration \begin{equation} \label{fibration} \Phi\: H_n\to
G(3,n) \end{equation} over the Grassmannian of 3-planes in $\PP^n$, with
fiber $H_3$.  It follows that $H_n$ is smooth and projective of
dimension $4n$.

There is a universal subscheme $\C\subset H_n\times \PP^n$.  For a
closed point $x\in H_n$, we denote by $C_x$ the corresponding cubic
curve, i.e., the fiber of $\C$ over $x$.  Also, let $\I_x\sub\OP{n}$ be
its ideal sheaf.  By the classification of the curves of $H_n$ (see
below), it is easy to see that
\begin{equation}
\label{H1}
	H^1(\PP^n,\I_x(d))=H^1(C_x,\OO_{C_x}(d))=0\quad
	\text{for all $d\ge1$ and for all
	$x\in H_n$.}
\end{equation}

For a subscheme $W\sub \PP^n$, denote by $H_W\sub H_n$ the closed
subscheme parameterizing twisted cubics contained in $W$.  There is a
natural scheme structure on $H_W$ as the intersection of $H_n$ with
the Hilbert scheme of $W$.  If $C_x\sub W$ is a Cohen-Macaulay twisted
cubic, then locally at $x\in H_n$, the scheme $H_W$ is simply the
Hilbert scheme of $W$.

Our goal is to compute the cycle class of $H_W$ in $A^*(H_n)$ in the
case that $W$ is a general complete intersection in $\PP^n$.  In
particular, we want its cardinality if it is finite, and its
Gromov-Witten invariants (see below) if it has positive dimension.

For each integer $d$ we define a sheaf $\E_{d}$ on $H_n$ by
\begin{equation}
	\label{defE}
	\E_{d}=p_{1*}(\OO_\C\* p_{2}^* \OP{n}(d)),
\end{equation}
where $p_1$ and $p_2$ are the two projections of $H_n\x \PP^n$.  If
$d\ge1$, then the vanishing of the first cohomology groups \eqref{H1}
implies by standard base change theory
\cite{AG} that $\E_{d}$ is locally free of rank $3d+1$, and moreover
that there are surjections $\rho\:H^0(\OP{n}(d))_{H_n}\to\E_{d}$ of
vector bundles on $H_n$.  In particular, for all $x\in H_n$, there is a
natural isomorphism
\begin{equation}
\E_{d}(x) @>\iso>> H^0(C_x,\OO_{C_x}(d)).
\label{basechange}
\end{equation}

A homogenous form $F\in H^0(\OP{n}(d))$ induces a global section
$\rho(F)$ of $\E_{d}$ over $H_n$, and the evaluation of this section at
a point $x$ corresponds under the identification (\ref{basechange}) to
the restriction of $F$ to the curve $C_x$.  Hence the zero locus of
$\rho(F)$ corresponds to the set of curves $C_x$ contained in the
hypersurface $V(F)$.

More generally, in the case of an intersection $V(F_1,\dots,F_p)$ in
$\PP^n$ of $p$ hypersurfaces, the section $(\rho(F_1),\dots,\rho(F_p))$
of $\E=\E_{d_1}\oplus\dots\oplus\E_{d_p}$ vanishes exactly on the points
corresponding to twisted cubics contained in $V(F_1,\dots,F_p)$.

\begin{prop} \label{propA}
Let $W\sub\PP^n$ be the complete intersection of $p$ general
hypersurfaces in $\PP^n$ of degrees $d_1,\ldots,d_p$ respectively.
Assume that $\sum_i(3d_i+1)=4n$.  Then the number of twisted cubic
curves contained in $W$ is finite and equals
$$\int_{H_n}c_{4n}(\E_{d_1}\oplus\cdots\oplus\E_{d_p}).$$
These cubics are all smooth.
\end{prop}

\begin{pf}

By the considerations above, the bundle
$\E=\E_{d_1}\oplus\cdots\oplus\E_{d_p}$ is a quotient bundle of the
trivial bundle $\dsum H^0(\OP{n}(d_i))_{H_n}$.  Hence Kleiman's Bertini
theorem \cite{Klei-1} implies that the zero scheme of the section
$(\rho(F_1),\ldots,\rho(F_p))$ is nonsingular and of codimension
$\rank(\E)$.  Since $\rank(\E)=\dim(H_n)$, the number of points is
finite and given by the top Chern class.
\end{pf}

\subsection{Gromov-Witten numbers}\label{B}
More generally, assume that $W$ is as in \propref{propA}, except that we
only assume an inequality $\sum_i(3d_i+1)\le4n$ instead of the equality.
The top Chern class of $\E$ still represents the locus $H_n(W)$ of
twisted cubics contained in $W$, although there are infinitely many of
them if the inequality is strict.  One may assign finite numbers to this
family by imposing Schubert conditions.  For this purpose say that a
{\em Schubert condition\/} on a curve is the condition that it intersect
a given linear subspace of $\PP^n$.  If the subspace has codimension
$c+1$, then the corresponding Schubert condition is of codimension $c$
(corresponding to the class $\gamma_c$ below).

\begin{defn}
Let $W$ be a general complete intersection in $\PP^n$ of
$p$ hypersurfaces of degrees $d_1,\dots,d_p$ respectively.  Let
$\lambda=(\lambda_1 \ge \lambda_2 \ge \dots \ge \lambda_m > 0)$ be a
partition of $4n - \sum_{i=1}^p(3d_i+1)$, and let $P_1,\dots,P_m$ be
general linear subspaces such that $\codim P_i = \lambda_i+1$.  The
number of twisted cubics on $W$ meeting all the $P_i$ is called the
$\lambda$-th {\em Gromov-Witten number} of the family of twisted
cubic curves on $W$, and is denoted by $N_3^\lambda(W)$.
\end{defn}

\begin{rem}
This is a slight variation on the definition used in \cite{Katz-2}, and
differs from that by a factor of 3 (resp.~9) for partitions with 2
(resp.~1) parts.  In \cite{Katz-2} only partitions of length at most three
are considered, as these numbers are the ones that have been predicted by
mirror symmetry computations (when $W$ is Calabi-Yau).  We have used the
term Gromov-Witten number rather than Gromov-Witten invariant, as the
latter term is now being used in a more sophisticated sense
\cite{Kont-Mani-1}.
\end{rem}

Let $h$ denote the hyperplane class of $\PP^n$ as well as its pullback to
$H_n\x\PP^n$, and let $[\C]\in A^*(H_n\x\PP^n)$ be the cycle class of the
universal curve $\C$.  If $P\sub\PP^n$ is a linear subspace of codimension
$c+1\ge2$, then $\C\cap H_n\x P$ projects birationally to its image under
the first projection, which is the locus of curves meeting $P$.  Hence the
class of the locus of curves meeting $P$ is ${p_1}_* (h^{c+1}[\C])$.  For
simplicity, we give this class a special notation:

\begin{notation}
For a natural number $c$, let
$\gamma_c = {p_1}_* (h^{c+1}[\C]) \in A^c(H_n)$.
\end{notation}

\begin{prop} \label{propB}
Let $W\sub\PP^n$ be the complete intersection of $p$ general
hypersurfaces in $\PP^n$ of degrees $d_1,\ldots,d_p$ respectively.
Assume that $\sum_i(3d_i+1)\le 4n$, and let $\lambda=(\lambda_1 \ge
\lambda_2 \ge \dots \ge \lambda_m > 0)$ be a partition of $4n -
\sum_{i=1}^p(3d_i+1)$.  Then
$$
N_3^\lambda(W)=
\int_{H_n} c_{{\text{top}}}(\E_{d_1}\oplus\cdots\oplus\E_{d_p})\cdot
\prod_{i=1}^m\gamma_{\lambda_i}.
$$
Furthermore, if $P_1,\dots,P_m$ are general linear subspaces such that
$\codim P_i = \lambda_i+1$, then the $N_3^\lambda(W)$ twisted cubics on
$W$ which meets each $P_i$ are all smooth.
\end{prop}

\begin{pf} Similar to the proof of \propref{propA}.
\end{pf}

For later use, we want to express the classes $\gamma_i$ in terms of
Chern classes of the bundles $\E_d$.

\begin{prop} \label{gammaformler}
Let $a_i=c_i(\E_1)$, $b_i=c_i(\E_2)$, $c_i=c_i(\E_3)$, and
$d_i=c_i(\E_4)$. Then we have the following formulas for the $\gamma_c$:
{\allowdisplaybreaks
\begin{align}
\gamma_0 =&\, 3\notag \\
\gamma_1 =&\, 5a_1-14b_1+13c_1-4d_1\notag\\
\gamma_2 =&\, 3a_1^2-9a_1b_1+9a_1c_1-3a_1d_1-3b_1^2+9b_1c_1\notag\\* &
-3b_1d_1-6c_1^2+3c_1d_1+a_2-3b_2+3c_2-d_2\notag\\
\gamma_3 =&\,
3a_1^3-9a_1^2b_1+9a_1^2c_1-3a_1^2d_1-3a_1b_1^2+9a_1b_1c_1\notag\\* &
-3a_1b_1d_1-6a_1c_1^2+3a_1c_1d_1-4a_1a_2-3a_1b_2+3a_1c_2-a_1d_2\notag\\* &+
14a_2b_1-13a_2c_1+4a_2d_1+3a_3\notag\\
\gamma_4 =&\, 3a_1^4-9a_1^3b_1+9a_1^3c_1-3a_1^3d_1-3a_1^2b_1^2\notag\\* &
+9a_1^2b_1c_1-3a_1^2b_1d_1-6a_1^2c_1^2+3a_1^2c_1d_1-7a_1^2a_2\notag\\* &
-3a_1^2b_2+3a_1^2c_2-a_1^2d_2+23a_1a_2b_1-22a_1a_2c_1+7a_1a_2d_1\notag\\* &
+3a_2b_1^2-9a_2b_1c_1+3a_2b_1d_1+6a_2c_1^2-3a_2c_1d_1+8a_1a_3\notag\\* &
-a_2^2+3a_2b_2-3a_2c_2+a_2d_2-14a_3b_1+13a_3c_1-4a_3d_1-3a_4\notag
\end{align}
}
\end{prop}

\begin{pf}
Let $\pi\: B=\PP(\E_1) \to H_n$.  The natural surjection
$\rho\:H^0(\OP{n}(1))_{H_n}\to\E_{1}$ induces a closed imbedding $B\sub
H_n\x\PP^n$ over $H_n$.  Over a closed point $x$ of $H_n$, the fiber of
$B$ is just the $\PP^3$ spanned by $C_x$.  So the universal curve $\C$
is actually a codimension 2 subscheme of $B$.  It follows by
the projection formula that
$$
\gamma_c = {\pi}_* (\tau^{c+1}[\C]_B) \in A^c(H_n),
$$
where $[\C]_B$ denotes the class of $\C$ in $A^2(B)$, and $\tau\in
A^1(B)$ is the first Chern class of the tautological quotient linebundle
on $B$.  The formulas of the proposition now follow by straightforward
computation (for example using \cite{schubert}) from the next lemma.
\end{pf}

\begin{lem} The class of $\C$ in $B$ is
 \begin{equation*}
 \begin{aligned}
 [\C]_B=&\,3\tau^2+(-4d_1+2a_1-14b_1+13c_1)\tau\\&
 +3c_1d_1+4a_2-3b_2+3c_2-d_2-2a_1^2\\&
 +5a_1b_1-4a_1c_1+a_1d_1-3b_1^2+9b_1c_1-3b_1d_1-6c_1^2
 \end{aligned}
 \end{equation*}
\end{lem}

\begin{pf}
Let $i\: \C\to B$ be the inclusion.  Then $[\C]_B$ equals
the degree 2 part of the Chern character of the $\OO_B$-module
$i_*\OO_\C(\ell\tau)$, for any integer $\ell$.  For $\ell=4$, there is a
canonical Beilinson type resolution of $i_*\OO_\C(4\tau)$:
\begin{equation}
	0\to\pi^*\E_1\*\Omega^3_{B/H_n}(3\tau)\to
	\pi^*\E_2\*\Omega^2_{B/H_n}(2\tau)\to
	\pi^*\E_3\*\Omega^1_{B/H_n}(\tau)\to
	\pi^*\E_4
	\label{beilinson}
\end{equation}
Using this it is a straightforward exercise (again using
\cite{schubert}) to compute the Chern character of $i_*\OO_\C$ in terms of
$\tau$ and the Chern classes of the $\E_n$.  \end{pf}

\subsection{Coarse classification of twisted cubics}\label{D}
We divide the curves $C_x$ for $x\in H_3$ into two groups, according to
whether they are Cohen-Macaulay or not.

A locally Cohen-Macaulay twisted cubic curve $C_x$ is also arithmetically
Cohen-Macaulay, and its ideal is given by the vanishing of the
$2\times2$ minors of a $3\x2$ matrix $\alpha$ with linear coefficients.
There is a resolution of $\OO_{C_x}$:
\begin{equation}
K^{\bullet}:\quad 0 @>>> F\*\OO_{\PP^3}(-3) @>\alpha>> E\*\OO_{\PP^3}(-2)
@>\wedge^2\alpha^t>> \OP3,
\end{equation}
where $F$ and $E$ are vector spaces of dimensions 2 and 3 respectively.
Intrinsically,
\begin{align}
    E=&H^0(\PP^3,\I_{x}(2))\sub H^0(\OP3(2))\label{EE}\\
    F=&\Ker(E\* H^0(\OP3(1)) @>\text{mult}>> H^0(\OP3(3))).\label{FF}
\end{align}

\begin{lem} \label{CMlem}
Let $C_x$ be Cohen-Macaulay, and let $E$ and $F$ be as above. Then
there is a functorial exact sequence:
\begin{equation}
0 @>>> \CC \to \End(F)\dsum\End(E) \to \Hom(F,E)\* H^0(\OP3(1))
\to T_{H_3}(x) \to 0
\end{equation}
\end{lem}

\begin{pf}
Recall the canonical identification
$T_{H_3}(x)=\Hom_{\PP^3}(\I_{x},\OO_{C_x})$.  The sequence now
follows from consideration of the total complex associated to the double
complex $\Hom_{\PP^3}(K^{\bullet},K^{\bullet})$.
\end{pf}

Next let us consider the curves $C_x$ for $x\in H_3$ which are not
Cohen-Macaulay.  By \cite{Pien-Schl}, these are projectively equivalent
to a curve with ideal generated by the net of quadrics
$x_0(x_0,x_1,x_2)$ plus a cubic form $q$, which can be taken to be of
the form $q=Ax_1^2+Bx_1x_2+Cx_2^2$, with $A$, $B$, and $C$ linear forms
in $\CC[x_1,x_2,x_3]$.  If we furthermore impose the conditions that $B$
is a scalar multiple of $x_3$, then the cubic $q$ is unique up to
scalar.  (See \cite{Elli-Stro-3}).

Let $Y\sub H_3$ be the locus of non-Cohen-Macaulay curves, and denote by
$I$ the 5-dimensional incidence correspondence $\{(p,H)\in
\PP^3\x{\PP^3}^* \mid p\in H \}$.  By the above, the quadratic part of
$\I_x$ for $x\in Y$ gives rise to a point of $I$.  This gives a morphism
\begin{equation}
	g\: Y\to I,
	\label{def:g}
\end{equation}
and again from the above it is clear that this makes $Y$ a
$\PP^6$-bundle over $I$.  Hence $Y$ is a divisor on $H_3$, and it is
clear how to compute the tangent spaces $T_Y(x)$.  To get hold of
$T_{H_3}(x)$, we need to identify the normal direction of $Y$ in $H_3$.

For this, let $C_x$ be the curve above, and consider the family $C_t$ of
Cohen-Macaulay curves given for $t\ne0$ by the matrix
\begin{equation}
\alpha_t =\begin{pmatrix} 0 & -x_0\\ x_0 & 0\\ -x_1 & x_2 \end{pmatrix}
         + t\begin{pmatrix} C & B\\  0 & A\\ 0 & 0 \end{pmatrix}
\end{equation}
Then
\begin{equation}
\det\left(\begin{array}{c|c} & x_1\\ \alpha_t & -x_2\\ & 0
\end{array}\right) = t(Ax_1^2+Bx_1x_2+Cx_2^2) = tq,
\end{equation}
which implies that $\lim_{t\to0}C_t=C_x$.  The tangent vector $\xi\in
T_{H_3}(x)=\Hom(\I_x,\OO_{C_x})$ corresponding to this one-parameter
family has this effect on the quadratic equations:
\begin{equation} \label{xi}
\xi(x_0^2)=Bx_0,\quad \xi(x_0x_1)=-Bx_1-Cx_2, \quad \xi(x_0x_2)=Ax_1.
\end{equation}
In particular, $\xi\ne0$. (This argument actually shows that the
blowup of the space of determinantal nets of quadrics
along the locus of degenerate nets maps isomorphically onto $H_3$,
cfr.~\cite{Elli-Pien-Stro}).

\subsection{The torus action}\label{E}
Consider the natural action of $\GL(n+1)$ on $\PP^n$.  It induces an
action on $H_n$ and on the bundles $\E_d$ for $d\ge1$.  Let
$T\sub\GL(n+1)$ be a maximal torus, and let $(x_0,\ldots,x_n)$ be
homogeneous coordinates on $\PP^n$ in which the action of $T$ is
diagonal.  A point $x\in H_n$ is fixed by $T$ if and only if the
corresponding curve $C_x$ is invariant under $T$, i.e., $t(C_x)=C_x$ for
any $t\in T$.  This is easily seen to be the case if and only if the
graded ideal of $C_x$ is generated by {\em monomials\/} in the $x_i$.
In particular, the fixpoints are isolated.

We will identify all the fixpoints $x\in H_n$, and for each of them we
will compute the representation on the tangent space $T_{H_n}(x)$.  The
tangent space of the Hilbert scheme is $\Hom(\I_x,\OO_x)$, but special
care must be taken at the points where $H_n$ meets another component
of the Hilbert scheme.  At these points, $T_{H_n}(x)$ is a proper
subspace of $\Hom(\I_x,\OO_x)$.

By the choice of the coordinates $(x_0,\ldots,x_n)$, there are
characters $\lambda_i$ on $T$ such that for any $t\in T$ we have
$t.x_i=\lambda^{}_i(t)x_i$.  The characters $\lambda^{}_i$ generate the
representation ring of $T$, i.e., if $W$ is any finite dimensional
representation of $T$ we may, by a slight abuse of notation, write
$W=\sum a_{p_0,\dots,p_n}\lambda_0^{p_0}\lambda_1^{p_1}
\cdots\lambda_n^{p_n}$, where the $p_i$ and $a_{p_0,\dots,p_n}$ are
integers.

Recall \eqref{fibration} the morphism $\Phi:H_n\to G(3,n)$ which maps a
point $x\in H_n$ to the 3-space spanned by the corresponding curve
$C_x$.  This morphism clearly is $\GL(n+1)$ equivariant, and its fibers
are all isomorphic to $H_3$.  If $C_x$ is invariant under $T$, then so
is its linear span.  Hence up to a permutation of the variables, we may
assume that it is given by the equations $x_4=\dots=x_n=0$, so that
$x_0,\dots,x_3$ are coordinates on the $\PP^3\sub\PP^n$ corresponding to
$\Phi(x)$.  The torus $T$ acts on $\PP^3$ through the four-dimensional
quotient torus $T_3$ of $T$ with character group spanned by
$\lambda_0,\dots,\lambda_3$.

The tangent space of $H_n$ at a fixpoint $x$ decomposes as a direct
sum
\begin{equation}
	\label{dsum}
	T_{H_n}(x)=T_{H_3}(x)\oplus T_{G(3,n)}(\Phi(x)),
\end{equation}
and it is well
known that
\begin{equation}\label{tangrass}
T_{G(3,n)}(\Phi(x))=
\Hom(H^0(\I_{\PP^3/\PP^n}(1)), H^0(\OP3(1)))=
\sum_{j=0}^3\sum_{i=4}^n
\lambda_j^{}\lambda_i\i.
\end{equation}
Hence we need to study the tangent space of $H_3=\Phi\i\Phi(x)$.

\begin{prop} Any fixpoint of $T_3$ in $H_3$ is projectively equivalent to
one of the following, where the first four are Cohen-Macaulay and the
last four are not:
\begin{alignat}{2}
(1)\quad& (x_0x_1,x_1x_2,x_2x_3) &\qquad \qquad
(5)\quad&(x_0^2,x_0x_1,x_0x_2,x_1x_2x_3) \notag\\
(2)\quad&(x_0x_1,x_1x_2,x_0x_2) &\qquad \qquad
(6)\quad&(x_0^2,x_0x_1,x_0x_2,x_1^{}x_2^2)\notag\\
(3)\quad&(x_0x_1,x_2^2,x_0x_2) &\qquad \qquad
(7)\quad&(x_0^2,x_0x_1,x_0x_2,x_2^2x_3^{})\notag\\
(4)\quad&(x_0^2,x_0x_1,x_1^2) &\qquad \qquad
(8)\quad&(x_0^2,x_0x_1,x_0x_2,x_2^3)\notag
\end{alignat}
\end{prop}

\begin{pf}
The action of $T_3$ on $\PP^3$ has the four coordinate points
as its fixpoints, and the only one-dimensional orbits are the six lines
of the coordinate tetrahedron.  Hence any curve invariant under $T_3$
must be supported on these lines.  If $C_x\in H_3$ is Cohen-Macaulay and
$T_3$-fixed, it is connected, has no embedded points and is not plane.
Hence there are only four possibilities: (1) the union of three distinct
coordinate lines, two of which are disjoint, (2) the union of three
concurrent coordinate lines, (3) a coordinate line doubled in a
coordinate plane plus a second line intersecting the first but not
contained in the plane, and finally (4) the full first-order
neighborhood of a coordinate line.

If $x\in Y^{T_3}$, then by the description of the curves in $Y$ we may
assume that the quadratic part of the ideal is $(x_0^2,x_0x_1,x_0x_2)$,
meaning that $C_x$ is a cubic plane curve in the plane $x_0=0$ which is
singular in $P=(0,0,0,1)$ plus an embedded point supported at $P$ but
not contained in the plane.  For the cubic we have these possibilities:
(5) the three coordinate lines, (6) one double coordinate line through
$P$ plus another simple line passing through $P$, (7) one double
coordinate line through $P$ plus another simple line not passing through
$P$, and (8) one coordinate line through $P$ tripled in the plane.
\end{pf}

\begin{rem}
There are several fixpoints of each isomorphism class, in
fact it is easy to verify by permuting the variables that in a given
$\PP^3$ the numbers of fixpoints of the types 1 through 8 are 12, 4, 24,
6, 12, 24, 24, 24, respectively.  This is consistent with the
fact that the (even) betti numbers of $H_3$ are 1, 2, 6, 10, 16, 19, 22,
19, 16, 10, 6, 2, 1, so that the Euler characteristic of $H_3$ is 130,
see \cite{Elli-Pien-Stro}.
\end{rem}

\begin{prop}\label{cm-akk}
Let $x$ be one of the fixpoints 1--4. Then the representation on the
tangent space $T_{H_3}(x)$ is given by
\begin{equation*}
T_{H_3}(x) = \Hom(F,E)\*(\l_0+\l_1+\l_2+\l_3)-\End(E)-\End(F)+1,
\end{equation*}
where the representations $E$ and $F$ are given in the following table:

\smallskip
\begin{center}
\begin{tabular} {|c|c|c|c|} \hline
Type\rommm & $\I_x$ & $E$ & $F$  \\ \hline
(1)\romm & $(x_0x_1,x_1x_2,x_2x_3)$ &
$\lambda^{}_0\lambda^{}_1+\lambda^{}_1\lambda^{}_2
+\lambda^{}_2\lambda^{}_3$ &
$\lambda^{}_0\lambda^{}_1\lambda^{}_2+\lambda^{}_1
\lambda^{}_2\lambda^{}_3$ \\ \hline
(2)\romm & $(x_0x_1,x_1x_2,x_0x_2)$ &
$\lambda^{}_0\lambda^{}_1+\lambda^{}_1\lambda^{}_2
+\lambda^{}_0\lambda^{}_2$ &
$2\lambda^{}_0\lambda^{}_1\lambda^{}_2$ \\ \hline
(3)\romm & $(x_0x_1,x_2^2,x_0x_2)$ &
$\lambda^{}_0\lambda^{}_1+\lambda_2^2+\lambda^{}_0
\lambda^{}_2$ &
$\lambda^{}_0\lambda^{}_1\lambda^{}_2+\lambda^{}_0
\lambda^{2}_2$ \\ \hline
(4)\romm & $(x_0^2,x_0x_1,x_1^2)$ &
$\lambda^{2}_0+\lambda^{}_0\lambda^{}_1+\lambda^{2}_1$ &
$\lambda^{}_0\lambda^{2}_1+\lambda^{2}_0\lambda^{}_1$  \\ \hline
\end{tabular}
\end{center}
\end{prop}

\begin{pf} Follows from \lemref{CMlem}, and the fact that $E$ and $F$
are equivariantly given by \eqref{EE} and \eqref{FF}.
\end{pf}

\begin{prop} \label{pr:akk}
Let $x$ be one of the fixpoints 5--8. Let $\mu$ be the character of
the minimal cubic generator, i.e., $\l_1\l_2\l_3$, $\l_1\l_2^2$,
$\l_2^2\l_3$, and $\l_2^3$, respectively, and let
\begin{align*}
	A=&\l_0\i(\l_1+\l_2+\l_3)+\l_3(\l_1\i+\l_2\i)\\
	B=&\l_1^3+\l_1^2\l_2+\l_1^2\l_3+\l_1\l_2^2+
	\l_1\l_2\l_3+\l_2^3+\l_2^2\l_3
\end{align*}
Then the tangent space of $H_3$ at $x$ is given by
\begin{equation*}
T_{H_3}(x) = A + \mu\i(B-\mu) + (\l_0\l_1\l_2)\i\mu
\end{equation*}
\end{prop}

\begin{pf}
Let $\beta=g(x) \in I$ be as in \secref{E}.  In fact, all
types 5--8 lie over the same fixpoint $\beta$.  The first term, $A$ in
the sum above, is easily seen to be the representation on the tangent
space $T_I(\beta)$.  Now $g\:Y\to I$ is a projective bundle, and the
fiber $g\i(\beta)$ is the projective space associated to the vector
space of cubic forms in $(x_1,x_2)^2\CC[x_1,x_2,x_3]$.  The
representation on this vector space is $B$, and the second term of the
formula of the proposition is the representation on $T_{g\i\beta}(x)$.
Thus the first two terms make up $T_Y(x)$.  The last term,
$(\l_0\l_1\l_2)\i\mu$, is the character on $N_{Y/H_3}(x)$.  This can be
seen from equations (\ref{xi}): by checking each case, one verifies that
the normal vector $\xi$ is semi-invariant with character
is $(\l_0\l_1\l_2)\i$ times the character of the cubic form $q$.
\end{pf}

\subsection{The computation}\label{F}
Let us briefly describe the actual computation, carried out using
``Maple'' \cite{Maple}, of the numbers in the introduction.  $H_n$ has a
natural torus action with isolated fixpoints.  By what we have done in
the last section, we can construct a list of all the fixpoints of $H_n$;
there are $130\binom{n+1}4 $ of these.  For each of them we compute the
corresponding tangent space representation, by \eqref{dsum} and
propositions \ref{cm-akk} and \ref{pr:akk}.

A consequence of the fact that all fixpoints are isolated is that none
of the tangent spaces contain the trivial one-dimensional
representation.  Choose a one-parameter subgroup $\psi\:\CC^* \to T$ of
the torus $T$, such that all the induced weights of the tangent space at
each fixpoint are non-zero.  This is possible since we only need to
avoid a finite number of hyperplanes in the lattice of one-parameter
subgroups of $T$.  For example, we may choose $\psi$ in such a way that
the weights of the homogeneous coordinates $x_0,\dots,x_n$ are
$1,w,w^2,\dots,w^n$ for a sufficiently large integer $w$.  In our
computations (for $n\le8$) we used instead weights taken from the sequence
4, 11, 17, 32, 55, 95, 160, 267, 441, but any choice that will not
produce a division by zero will do.

Since all the tangent weights of the $\CC^*$ action on $H_n$ so obtained
are non-trivial, it follows that this action has the same fixpoints as
the action of $T$, hence a finite number.

By \propref{propB}, we need to evaluate the class
$$\delta=c_{{\text{top}}}(\E_{d_1}\oplus\cdots\oplus\E_{d_p})\cdot
\prod_{i=1}^m\gamma_{\lambda_i}\in A^{4n}(H_n).$$

Note that the isomorphism \eqref{basechange}
is equivariant.  Clearly $H^0(\OO_{C_x}(d))$ is spanned by all monomials
of degree $d$ not divisible by any monomial generator of $I_x$.  Thus we
know all the representations $\E_d(x)$ for all fixpoints $x\in H_n$.

By \propref{gammaformler}, $\delta$ is a polynomial
$p(\dots,c_k(\E_d),\dots)$ in the Chern classes of the equivariant
vector bundles $\E_d$.  To find a function $f$ on the fixpoint set which
represents $\delta$, simply replace each occurance $c_k(\E_d)$ by the
localized equivariant Chern class $\sigma_k(\E_d,-)$, i.e.,  put
$f=p(\dots,\sigma_k(\E_d,-),\dots)$.  Then the class is evaluated by
the formula in \thmref{bott1} (2).

\section{The Hilbert scheme of points in the plane}\label{hilbkap}

Let $V$ be a three-dimensional vector space over $\CC$ and let
$\PP(V)$ be the associated projective plane of rank-1 quotients of $V$.
Denote by $H_r=\Hilb^r_{\PP(V)}$ the Hilbert scheme parameterizing
length-$r$ subschemes of $\PP(V)$. There is a universal subscheme
$\Z\sub H_r\x\PP(V)$. We will use similar notational conventions
as in \secref{A}: for example, if $x\in{H_r}$, the corresponding subscheme
of $\PP^2$ is denoted $Z_x$, its ideal sheaf $\I_x$ etc.
As in \eqref{defE}, let for any integer $n$
\begin{equation}
	\label{def2E}
	\E_{n}=p_{1*}(\OO_\Z\* p_{2}^* \OO_{\PP(V)}(n)),
\end{equation}
where $p_1$ and $p_2$ are
the two projections of $H_r\x\PP(V)$.
Since $\Z$ is finite over ${H_r}$, all basechange maps
\begin{equation}
\E_{n}(x) @>\iso>> H^0(Z_x,\OO_{Z_x}(n)).
\label{basechange2}
\end{equation}
are isomorphisms. In particular, $\E_n$ is a rank-$r$ vector bundle
on ${H_r}$. Denote by $\L$ the linebundle
\begin{equation}\label{defL}
\L = \ytre^r\E_0\*\ytre^r\E_{-1}\v.
\end{equation}
Then $\L$ corresponds to the divisor on $H_r$ corresponding to
subschemes $Z$ meeting a given line.  We are going to compute integrals
of the form
\begin{equation}\label{generalintegral}
\int_{{H_r}} s_{2r}(\E_n\*\L^{\*m})
\end{equation}
for small values of $r$.  Afterwards we will give interpretations of
some of these numbers in terms of degrees of power sum and Darboux loci
in the system $\PP(S_nV)$ of plane curves of degree $n$ in the dual
projective plane $\PP(V\v)$.

As usual, we start by identifing all the fixpoints and tangent space
representations for a suitable torus action.  This has been carried out
in more detail in \cite{Elli-Stro-1}, the following simpler presentation
is sufficient for the present purpose.

As in \secref{E}, let $T\sub \GL(V)$ be a maximal torus and let
$x_0,x_1,x_2$ be a basis of $V$ diagonalizing $T$ under the natural
linear action.  The eigenvalue of $x_i$ is a character $\l_i$ of $T$.
We identify characters with one-dimensional representations, hence the
representation ring of $T$ with the ring of Laurent polynomials in
$\l_0,\l_1,\l_2$.  For example, the natural representation on the vector
space $V\v$ can be written $\l_0\i+\l_1\i+\l_2\i$.

Fixpoints of ${H_r}$ can be described in terms of partitions,
i.e., integer sequences $b=\{b_r\}_{r\ge0}$ weakly decreasing to zero.
Let $|b|=\sum_{r\ge0}b_r$. The {\em diagram\/} of a partition $b$ is
the set $D(b)=\{(r,s)\in\ZZ_{\ge0}^2\mid s< b_r\}$ of cardinality $|b|$.
 A {\em
tripartition\/} is a triple $B=(b^{(0)},b^{(1)},b^{(2)})$ of
partitions; put $|B|=\sum_i |b^{(i)}|$; the number being partitioned.
The {\em $n$-th diagram\/} $D_n(B)$ of a tripartition
$B=(b^{(0)},b^{(1)},b^{(2)})$ is defined for
$n\ge|B|$ as follows: Letting the index $i$ be counted modulo 3, we put
$$
D_n^i(B)=\{ (n_0,n_1,n_2)\in\ZZ^3 \mid n_0+n_1+n_2=n
\text{ and } (n_{i+1},n_{i+2})\in
D(b^{(i)})\}
$$
and $D_n(B) = D_n^0(B)\cup D_n^1(B)\cup D_n^2(B)$.  Intuitively, the
diagram $D_n(B)$ lives in an equilateral triangle with corners $(n,0,0)$,
$(0,n,0)$, and $(0,0,n)$, and originating from the $i$-th corner there is a
(slanted) copy of $D(b^{(i)})$.  When $n\ge|B|$, these don't overlap.  As
$n$ grows, the shape and size of the three parts of $D_n(B)$ stay the same,
whereas the area separating them grows.  We may also define $D_n(B)$ for
integers $n<r$ by the same formula as above, but where the union is taken
in the sense of multisets, i.e., some elements might have multiplicities 2
or even 3.  For $n<r$ the diagram $D_n(B)$ may also stick out of the
triangle referred to above.

A fixpoint $x\in{H_r}$ corresponds to a length-$r$ subscheme $Z_x\sub\PP(V)$
defined by a monomial ideal.  Fix an integer $n\ge r$ and consider the
set
$$
D_n(Z_x)=\{(n_0,n_1,n_2)\in \ZZ_{\ge0}^3 \mid \sum_i n_i=n \text{ and }
\prod x_i^{n_i} \notin H^0(\PP(V),\I_x(n))\}
$$
This set is the $n$-th diagram of a tripartition $B$ of $r$, the three
constituent partitions corresponding to the parts of $Z_x$ supported in the
three fixpoints of $\PP(V)$.  Conversely, starting with a tripartition of
$r$, we may obviously construct a monomial ideal of colength $r$ in
such a way that we get an inverse of the construction above. Hence there
is a natural bijection between $H_r^T$ and the set of tripartitions of
$r$.

In terms of representations, it follows from the above that for a
fixpoint $x$ corresponding to the tripartition
$B=(b^{(0)},b^{(1)},b^{(2)})$, we have
\begin{equation}\label{reprEn}
\E_n(x)=H^0(\OO_{Z_x}(n)) = \sum_{(n_0,n_1,n_2)\in D_n(B)} \prod \l_i^{n_i}.
\end{equation}
For $n<r$, the summation index needs to be interpreted as running
through the multiset $D_n(B)$.  The representation on $\L$ in the same
fixpoint is
\begin{equation}\label{reprL}
\L(x) = \prod \l_i^{|b^{(i)}|}.
\end{equation}
For $n\ge r$, we also have the following formula for
$I_n:=H^0(\PP(V),\I_x(n))$ in the representation ring:
$$
I_n =S_nV - H^0(\OO_{Z_x}).
$$

To compute the tangent space representation, we use a trick that is often
useful even in higher dimensions: functorial free resolutions.  The tangent
space of ${H_r}$ in $x$ is canonically isomorphic to $\Ext^1(\I_x,\I_x)$.
Fix an integer $n\ge r+2$.  Then there is a canonical resolution of locally
free $\OO_{\PP(V)}$-modules
$$
 K_{\bullet}:\quad 0\to K_2 \to K_1 \to K_0
$$
of
$\I_x(n)$, where $K_p=\Omega_{\PP(V)}^p(p)\*I_{n-p}$.  As in
\eqref{beilinson}, this is a special case of Beilinson's spectral sequence.
$T$ acts on $K_\bullet$.  Let $S^\bullet$ be the total complex
associated to the double complex $\Hom_{\PP^2}(K_\bullet,K_\bullet)$.
Then the $i$-th cohomology group of $S^\bullet$ is
$\Ext^i(\I_x(n),\I_x(n))=\Ext^i(\I_x,\I_x)$ for $i=0,1,2$
\cite[Lemma 2.2]{Elli-Stro-4}.
For $i=0$ this is $\CC$ (with trivial action) and for
$i=2$ it is zero.  Using the canonical identifications
$\Hom(\Omega^p(p),\Omega^q(q)) = \ytre^{p-q}V\v$, we end up with the
following formula for the tangent space representation in terms of the
data of the tripartition $B$:

{\allowdisplaybreaks
\begin{align}\label{reprT}
T_{{H_r}}(x)  = \,& 1-(\sum_{i=0}^2(-1)^i\Ext^i(\I_x,\I_x))\\
= \,&1 -(S^0-S^1+S^2) \notag\\
=\,&1-(\Hom(I_n,I_{n})+\Hom(I_{n-1},I_{n-1})+\Hom(I_{n-2},I_{n-2}))\notag\\*
&+ (\l_0\i+\l_1\i+\l_2\i)
(\Hom(I_{n-1},I_n)+\Hom(I_{n-2},I_{n-1}))\notag\\*
&-(\l_0\i\l_1\i+\l_1\i\l_2\i+\l_2\i\l_0\i)\Hom(I_{n-2},I_n)\notag
\end{align}
}

Here are the computational results on the Hilbert scheme which will
be used in the following applications:
\begin{prop} \label{computation1} Let $\E_n$ be as in \eqref{def2E}.
For $2\le r\le 8$, we have
$$
\int_{{H_r}}s_{2r} (\E_n) = s_r(n),
$$
where $s_r(n)$ are as in \thmref{main4}.
\end{prop}

\begin{prop} \label{computation2} Let $\E_{-1}$ and $\L$ be as in
\eqref{def2E} and \eqref{defL}.
For $r=2,3,4,5,6,7,8,9,10$, we have
$$
\int_{{H_r}}s_{2r}( \E_{-1}\*\L) = 0,0,0,0,2540, 583020, 99951390,
16059395240, 2598958192572.
$$
\end{prop}

\begin{pf} For both propositions, apply Bott's formula. The one-parameter
subgroup of $T$ such that the $\l_i$ have weights 0,1,19 will work.  The
contribution at each fixpoint is given by
\eqref{reprEn}, \eqref{reprL}, and \eqref{reprT}. Generate all
fixpoints and perform the summation using e.g. \cite{Maple}.
\end{pf}

\begin{rem} There is no difficulty in principle to evaluate
\eqref{generalintegral} directly with symbolic values of both $n$ and
$m$, for given values of $r$. For example, for $r=3$, the result is
$(n^{6}+24\,n^{5}m+252\,n^{4}m^{2}+1344\,n^{3}m^{3}+3780\,
n^{2}m^{4}+5040\,nm^{5}+
2520\,m^{6}-30\,n^{4}-432\,n^{3}m-2520\,n^{2}m
^{2}-6048\,nm^{3}-5040\,m^{4}+45\,n^{3}+ 504\,n^{2}m+2268\,nm^{2}+3024
\,m^{3}+206\,n^{2}+1200\,nm+1512\,m^{2}-576\,n- 1728\,m+384)/6$. However,
with given computer resources, one gets further the fewer variables one
needs. On a midrange workstation, we could do this integral up to $r=5$.
\end{rem}

\begin{rem} Tyurin and Tikhomirov \cite{Tyur-Tikh} and Le Potier have
shown that \propref{computation2} implies that the Donaldson
polynomial $q_{17} (\PP^2) = 2540 $.  It may also be deduced from the
proposition that $q_{21}(\PP^2) = 233208$, see \cite{Tyur-Tikh} or our
forthcoming joint paper with J. Le Potier.
\end{rem}

\subsection{Power sum varieties of plane curves}
Closed points of $\PP(S_nV)$ correspond naturally to curves of degree
$n$ in the dual projective plane $\PP(V\v)$.  In particular, points of
$\PP(V)$ correspond to lines in $\PP(V\v)$, so $H_r=\Hilb^r_{\PP(V)}$
is a compactification of the set of unordered $r$-tuples of linear
forms modulo scalars.

Let $r$ and $n$ be given integers.  Let $U(r,n)$ be the set of pairs
$(\{L_1,\dots,L_r\},C)$ where the $L_i\sub\PP(V\v)$ are lines in
general position, and $C$ is a curve with equation of the form
$\sum_{i=1}^r a_i\ell_i^n \in S_nV$, where $\ell_i\in V\v$ is an
equation of $L_i$.  Then the power sum variety $PS(r,n)$ is the
closure of the image of $U(r,n)$ in $\PP(S_nV)$ under the projection
onto the last factor.  To compute the degree of the image times the
degree $p(r,n)$ of the map $U(r,n)\to PS(r,n)$, we need to find a
workable compactification of $U(r,n)$.

Recall from \eqref{def2E} the rank-$r$ vector bundle $\E_n$ on
${H_r}$.  It comes naturally with a morphism $S_nV_{H_r}\to \E_n$,
which is surjective if $n\ge r-1$.  Now consider the projective bundle
over $H_r$:
\begin{equation}
\PP(\E_n) \sub \PP(S_nV) \x {H_r}.
\end{equation}
It is easy to verify that $\PP(\E_n)$ contains $U(r,n)$ as an open
subset.  It follows that $p(r,n)$ times the degree of $PS(r,n)$ is given
by the self-intersection of the pullback to $\PP(\E_n)$ of
$\OO_{\PP(S_nV)}(1)$. This is
$\int_{\PP(\E_n)} c_1(\OO_{\PP(\E_n)}(1))^{3r-1}$, and pushing it
forward to the Hilbert scheme, we get almost by definition,
$\int_{{H_r}} s_{2r}(\E_n)$
\cite[Ch.~ 3]{IT}.
Together with \propref{computation1}, this proves \thmref{main4}.

\subsection{Darboux curves} These curves are also defined in terms of
linear forms, and we may take $H_{n+1}$ as a parameter space for the
variety of complete $(n+1)$-gons.  For a length-$(n+1)$ subscheme
$Z\sub\PP(V)$, put $E=H^0(\OO_{Z}(-1))=H^1(\I_Z(-1))$ and $F=H^1(\I_Z)$.
The multiplication map
 $
 V\* E \to F
 $
 gives rise to a bundle map over the dual plane $\PP(V\v)$:
$$
m\: E_{\PP(V\v)}(-1) \to F_{\PP(V\v)}.
$$
If $Z$ consists of $n+1$ points in general position,
the degeneration locus $D(Z)$ corresponds to the set of bisecant
lines to $Z$, i.e., the singular locus of the associated $(n+1)$-gon.
The Eagon-Northcott resolution of $D(Z)$ gives the
following short exact sequence:
$$
0 \to F_{\PP(V\v)}\v(-1) \to E_{\PP(V\v)}\v \to L\*_\CC\I_{D(Z)}(n) \to 0,
$$
showing that there is a natural surjection
$$
S_nV\v \to H^0(\PP(V\v),\I_{D(Z)}(n))\v \iso H^0(Z,\OO_Z(-1))\*L.
$$
Here $L$ is
the onedimensional vector space $\det(F)\*\det(E)\i$.

Globalizing this construction over $H_{n+1}$ gives a natural map
$$
p\: S_nV_{H_{n+1}} \to \E_{-1}\*\L
$$
such that the closure of the image of the induced rational map
$\PP(\E_{-1}\*\L)\to\PP(S_nV)$ is the Darboux locus $D(n)$.
By the lemma below, this rational map is actually a morphism. Thus
we may argue as in the power sum case and find that $p(n)$ times the
degree of $D(n)$ is $\int_{\PP(\E_{-1}\*\L)} c_1(\OO(1))^{3n+2}=
\int_{H_{n+1}} s_{2n+2}(\E_{-1}\*\L)$. This together with
\propref{computation2} implies \thmref{main5}.

\begin{lem}
The  bundle map $p\:S_nV_{H_{n+1}} \to \E_{-1}\*\L$ over
$H_{n+1}$ is surjective.
\end{lem}

\begin{pf} Assume the contrary. Since the support of the cokernel is
closed and $\GL(V)$-invariant,
there exists a subscheme $Z$ in $\Supp\Coker(p)$ which is supported in
one point.  Without loss of generality we may assume that $Z$ is
supported in the point $x_1=x_2=0$.

Let $E$ and $F$ be as above, and let $K\sub E$ be a subspace of
codimension 1.  For a linear form $\ell\in V$, let $m_{\ell}\:E\to F$
be multiplication by $\ell$.  The assumption that $p$ is not
surjective means that $K$ can be chosen such that the determinant of
the restriction of $m_{\ell}$ to $K$ is 0 for all $\ell\in V$.

Multiplication by $x_0$ induces an isomorphism $E=H^0(\OO_Z(-1)) \iso
H^0(\OO_Z) \iso \CC[x,y]/I_Z$, where $x=x_1/x_0$ and $y=x_2/x_0$.  Under
these identifications, if $\ell=1-ax-by$ is the image of a general
linear form, the kernel of $m_{\ell}$ is generated by $1/\ell$.

Consider the set $S$ consisting of all such elements $\ell\i$, with
$a,b\in\CC$.  The series expansion $\ell\i=1+(ax+by)+(ax+by)^2+\cdots$
shows that $S$ generates $\CC[x,y]/(x,y)^m$ as a vector space for all
$m\ge0$.  Indeed, a hyperplane $W_m\sub\CC[x,y]/(x,y)^m$ containing the
image of $S$ would, by induction on $m$, dominate
$\CC[x,y]/(x,y)^{m-1}$.  Hence $W_m$ could not contain
$(x,y)^{m-1}/(x,y)^m$.  But the image of $S$ in $(x,y)^{m-1}/(x,y)^m$ is
the cone over a rational normal curve of degree $m-1$, hence spans this
space.

Since $(x,y)^m\sub I_Z$ for $m$ large, it follows that $S$ generates
$\OO_Z$ and hence $E$ as a vector space.  Now for an $\ell$ such that
$\ell\i\notin K$, the restriction of $m_\ell$ to $K$ will be an
isomorphism.  This leads to the desired contradiction.
\end{pf}

\section{Discussion} How general is the strategy of using Bott's
formula in enumerative geometry, as outlined in these examples? The
first necessary condition is probably a torus action, although Bott's
formula is valid in a more general situation: a vector field with
zeros and vector bundles acted on by the vector field.  It seems to
us, though, that the cases where one stands a chance of analysing the
local behaviour of bundles near all zeroes of such a field are those
where the both the vector field and its action on the vector bundles
are ``natural'' in some sense.  If there are parameter spaces with
natural flows on them, not necessarily coming from torus actions,
presumably Bott's formula could be useful.

It is not necessary that the fixpoints be isolated in order for the
method to give results.  Kontsevich's work \cite{Kont-1} is a
significant example of this.  Another natural candidate for Bott's
formula is the moduli spaces of semistable torsionfree sheaves on
$\PP^2$.  These admit torus actions, but not all fixpoints are
isolated.  One can still control the structure of the fixpoint
components, however.  This may hopefully be used for computing
Donaldson polynomials of the projective plane, at least in some cases.

A more serious obstacle to the use of Bott's formula is the presence
of singularities in the parameter space.  For example, all components
of the Hilbert scheme of $\PP^n$ admit torus actions with isolated
fixpoints, but they are almost all singular.  The main non-trivial
exceptions are actually the ones treated in the present paper.  On the
other hand, singularities present inherent problems for most
enumerative approaches, especially if a natural resolution of
singularities is hard to find.

For all examples in this paper, one needs a computer to actually perform
the tedious computations.  We mentioned already that the number of
fixpoints in the case of twisted cubics in $\PP^n$ is $130\binom{n+1}4$.
For the Hilbert scheme of length-8 subschemes of $\PP^2$, the
corresponding number is 810.  From the point of view of computer
efficiency, there are some advantages to the use of Bott's formula in
contrast to trying to work symbolically with generators and relations in
cohomology, as for example in \cite{Elli-Stro-3} or \cite{schubert}.
First of all, the method works even if we don't know all relations, as
is the case for the Hilbert scheme of the plane, for example.  But the
main advantage is perhaps that Bott's formula is not excessively hungry
for computer resources.  It is often straightforward to make a loop over
all the fixpoints.  The computation for each fixpoint is fairly simple,
and the result to remember is just a rational number.  This means that
the computer memory needed does not grow much with the number of
fixpoints, although of course the number of CPU cycles does.  For
example, most of the numbers of elliptic quartics were computed on a
modest notebook computer, running for several days.

Finally, Bott's formula has a nice error-detecting feature, which is
an important practical consideration: If your computer program
actually produces an integer rather than just a rational number,
chances are good that the program is correct!


\end{document}